\documentclass[twocolumn]{svjour3}  
\smartqed  

\usepackage[dvipdfmx]{color}
\usepackage{graphicx}
\usepackage[T1]{fontenc}
\usepackage{mathptmx}  
\usepackage[scaled]{helvet}  
\usepackage{courier}
\usepackage[subrefformat=parens]{subcaption}

\makeatletter
\let\cl@chapter\undefined
\makeatletter
\usepackage{amsmath,amssymb}   
\usepackage{mathtools} 
\usepackage{cleveref}  
\Crefname{equation}{Eq.}{Eqs.}%
\Crefname{figure}{Fig.}{Figs.}%
\usepackage{amsbsy}
\usepackage{bm}
\usepackage{algorithmic}
\usepackage{algorithm}

\def\makeheadbox{\relax}
\begin{document}

\title{Stochastic Assessment of Acceleration Probability Density Function for Parametric Rolling Using Moment Method
}

\author{Yuuki Maruyama      \and Atsuo Maki \and
        Leo Dostal          \and Naoya Umeda
}

\institute{Yuuki Maruyama \and Atsuo Maki \and Naoya Umeda \at
              Department of Naval Architecture and Ocean Engineering, Graduate School of Engineering, Osaka University, 2-1 Yamadaoka, Suita, Osaka, Japan \\
              \email{yuuki\_maruyama@naoe.eng.osaka-u.ac.jp} 
           \and
           Leo Dostal \at
           Institute of Mechanics and Ocean Engineering, Hamburg University of Technology, 21043 Hamburg, Germany
}

\date{Received: date / Accepted: date}

\maketitle

\begin{abstract}
Container ships encounter large roll angles and high acceleration, and container loss remains a problem. This study proposes a method for calculating the probability density function~(PDF) of roll angular and cargo lateral accelerations. First, the moment values of these accelerations are derived using the linearity of expectation and the validity of this method is examined. Second, the PDF shapes of these accelerations are proposed and their coefficients are determined using the obtained moment values. Our proposed method can be used to derive the PDFs of roll angular and cargo lateral accelerations.

\end{abstract}

%
\section{Introduction}
%
%

The second-generation intact stability criteria were developed by the International Maritime Organization~(IMO)~\cite{IMO2020}. 
Container ships encounter large roll angles and high acceleration, and container loss remains a problem. The failure modes relevant to container loss accidents are parametric rolling and excessive acceleration. 
The new criteria in direct stability assessment require using the direct counting method. 
However, this simulation is time-consuming. Thus, a theoretical estimation method is significantly desirable as a preliminary tool.

Because parametric rolling is a nonlinear phenomenon, estimating it simply and mathematically is difficult. One solution is to use probability theory. Many studies~\cite{Roberts1982para}\cite{Dostal2012}\cite{Maruyama2021}\cite{Maruyama2022moment} have been conducted on the probability density function~(PDF) of roll angle and amplitude for parametric rolling. For instance, Maki et al.~\cite{Maki2021} proposed the “PDF Line Integral method~(PLIM)” and obtained the PDF of the acceleration’s main component in beam seas. Maruyama et al.~\cite{Maruyama2022accel} investigated the PDF of acceleration for parametric rolling in longitudinal waves. Here, roll angular acceleration was divided into two components, and each PDF was derived using PLIM. Both studies could not simply obtain the PDF of the roll angular acceleration. Furthermore, little research has been conducted on the PDF of roll angular and cargo lateral accelerations for parametric rolling. 

In this study, we propose a method~[Sec.~\ref{sec:Mathematical Model}] for calculating the moment values of roll angular and cargo lateral accelerations. Furthermore, the PDFs of these accelerations are derived using their moment values~[Sec.~\ref{sec:sec4sub1}]. Here, the PDF shape is suggested, and their coefficients are optimized~[Sec.~\ref{sec:sec4sub2}].

%
\section{Subject Ship and Sea Condition}
\label{sec:Subject Ship and Sea Condition}
\begin{figure}[h]
    \centering
    \includegraphics[scale=1]{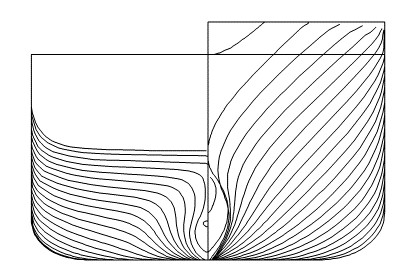}
    \caption{Body Plan (C11)}
    \label{fig:bodyplanC11}
\end{figure}

\begin{figure}[h]
    \centering
    \includegraphics[scale=0.8]{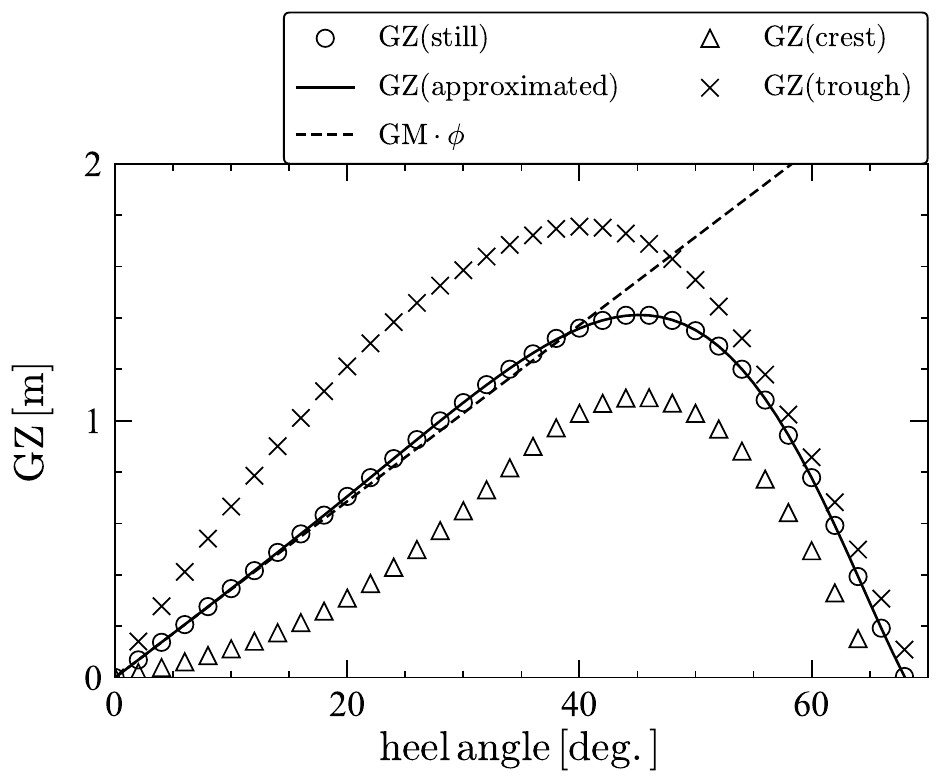}
    \caption{GZ curve in still water and for wave crest and wave trough conditions (C11); this figure was adapted from Maruyama et al.~\cite{Maruyama2022accel} }
    \label{fig:GZcurve_stillwater_C11}
\end{figure}

In this study, the results of the calculation for Froude number Fn = 0.00 in head seas are compared. Please note that in the considered case the choice of this Froude number leads to most severe parametric rolling conditions. The subject ship is a post-Panamax container ship of the C11 class~\cite{Hashimoto2010}. This ship's body plan and major characteristics are shown in Fig.~\ref{fig:bodyplanC11} and Table~\ref{tab:principal_C11}, respectively.
In this figure, the circle represents the GZ curve in calm water. Additionally, the triangle and cross marker illustrate the GZ curves obtained from the hydrostatic calculation in case the wave crest or trough is located at amidships.
The wave is a regular wave with the same wavelength as the ship, in this case. 
These restoring arms is calculated based on the Froude-Krylov hypothesis~\cite{Hamamoto1991}. 
The actual GZ is approximated using a 9th order polynomial to obtain a suitable GZ curve. The C11 container ship has a linear GZ curve up to a roll angle of approximately 40 degrees, as shown in Fig.~\ref{fig:GZcurve_stillwater_C11}. Furthermore, the Ikeda’s simplified method~\cite{Kawahara2012} is used to estimate the roll damping coefficients. Therefore, the damping coefficients are $\beta_1^{}=3.64\times10_{}^{-3}$ and $\beta_3^{}=4.25$.
\renewcommand{\arraystretch}{1.5}
\begin{table}[h]
  \begin{center}
    \caption{Principal particulars of the subject ship at full scale}
    \begin{tabular}{ll}
    \hline
     {\it Items} & {\it C11}
     \\
     \hline
      Length:\,$L_{pp}$ & $262.0 \,[\mathrm{m}]$
      \\
      Breadth:\,$B$ & $40.0 \,[\mathrm{m}]$
      \\
      Depth:\,$D$ & $24.45 \,[\mathrm{m}]$
      \\
      Draught:\,$d$ & $11.5 \,[\mathrm{m}]$
      \\
      Block coefficient:\,$C_b$ & $0.562$
      \\
      Metacentric height:\,$\mathrm{GM}$ & $1.965 \,[\mathrm{m}]$
      \\
      Natural roll period:\,$T_{\phi}$ & $25.1 \,[\mathrm{s}]$
      \\
      Bilge keel length ratio:\,$L_{BK}/L_{pp}$ & $0.292$
      \\
      Bilge keel breadth ratio:\,$B_{BK}/B$ & $0.0100$
      \\
    \hline
    \end{tabular}
    \label{tab:principal_C11}
  \end{center}
\end{table}
\renewcommand{\arraystretch}{1.0}


In this study, the ITTC spectrum is used to approximate the ocean wave spectrum, which is given by
\begin{equation}
    S_{\mathrm{w}}(\,\omega\,) = \dfrac{173\,H_{1/3}^{2}}{T_{01}^{4}\,\omega^{5}}\, \exp\left(\,-\,\dfrac{691}{T_{01}^{4}\,\omega^{4}}\,\right).
\end{equation}
In this study the wave mean period $T_{01}^{} = 9.99$[s] and the significant wave height $H_{1/3}^{} = 5.0$[m] are used. In this study, the GM variation is calculated by using the concept of Grim’s effective wave \cite{Grim1961} to irregular waves. The transfer function $H_{\zeta}^{}$, which must obtain the time series of the effective wave $A_{\text{w}}^{}(t)$, can be represented as shown in the following equation\cite{Umeda1991}:
\begin{equation}
    \begin{split}
        \displaystyle
        & H_{\zeta}^{}(\,\omega, \chi\,) = H_{\zeta c}^{}(\,\omega, \chi\,) + i H_{\zeta s}^{}(\,\omega, \chi\,)
        \\\\
        & \left\{\begin{split}
        \displaystyle
        & H_{\zeta c}^{}(\,\omega, \chi\,) = \dfrac{\dfrac{\omega^{2}\,L}{g}\,\cos\chi\, \sin\left(\,\dfrac{\omega^{2}\,L}{2\,g}\,\cos\chi\,\right)}{\pi^{2} - \left(\,\dfrac{\omega^{2}\,L}{2\,g}\,\cos\chi\,\right)^{2}}
        \\\\
        \displaystyle
        & H_{\zeta s}^{}(\,\omega, \chi\,) = 0
        \\\\
        \end{split}\right.
    \end{split}
\end{equation}
Here, $\omega$ denotes the wave angular frequency, $\chi$ denotes the heading angle from wave direction, $g$ denotes gravitational acceleration, and $L$ denotes the wave length.
Using this transfer function, the spectrum of the effective wave can be obtained as follows:
    \begin{equation}
        \label{spectrum_effective_wave_theory}
        S_{\mathrm{eff}}^{}(\,\omega\,) = \left|\,H_{\zeta}^{}(\,\omega\,)\,\right|^{2} S_{\mathrm{w}}^{}(\,\omega\,)
    \end{equation}
To prevent time series repetition, a technique \cite{Shuku1979} is used to divide the spectrum into elements with equal wave energies.

Moreover, the Monte Carlo simulation~(MCS) is computed using the fourth-order Runge–Kutta method. A time step of 0.02[s] is specified. The initial conditions are 5[deg.] roll angle and 0[deg./s] roll velocity. The number of realizations is 1000, and each simulation lasts 1 h. 

%
\section{Mathematical Model}
\label{sec:Mathematical Model}
\begin{figure}[h]
    \centering
    \includegraphics[scale=1]{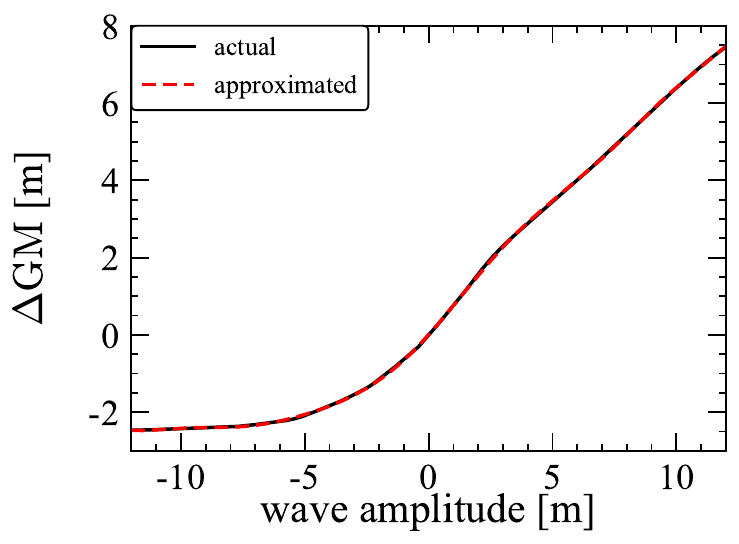}
    \caption{Relationship between $\Delta$GM and wave amplitude at amidships, subject ship : C11. This figure was adapted from Maruyama et al.~\cite{Maruyama2022moment} }
    \label{fig:HBL_vs_GM}
\end{figure}

In this study, Eq.~\eqref{1DoF_roll_equation} represents the equation for roll motion in irregular waves. 
Here, the additive wave excitation $f(t)$ is absent ($\,f(\,t\,)=0\,$), because longitudinal wave conditions only lead to parametric excitation. 
\begin{equation}
    \label{1DoF_roll_equation}
    \ddot{\phi} + \beta_1^{}\,\dot{\phi} + \beta_3^{}\,\dot{\phi}_{}^3 + \sum_{n=1}^{5} \alpha_{2n-1}^{}\,\phi_{}^{2n-1} + P(\,t\,)\,\phi = f(\,t\,)
\end{equation}
Here, the roll angle, roll velocity, and roll angular acceleration are denoted by $\phi$, $\dot{\phi}$, and $\ddot{\phi}$, respectively. The parameter $\beta_{1}^{}$ is linear, and $\beta_{3}^{}$ is the cubic damping coefficient, divided by $I_{xx}^{}$, where $I_{xx}^{}$ denotes the moment of inertia in roll (including the corresponding added moment of inertia), and $\alpha_{i}^{}\,(i=1,3,5,7, 9)$ denotes the coefficient of the i-th component of the polynomial fitted GZ curve. Additionally, $P(\,t\,)$ denotes a time series of parametric excitation as follows:
\begin{equation}
\label{non-memory transformation}
    \displaystyle
    P(\,t\,) = \frac{\omega_{0}^{2}}{GM} \sum_{n=1}^{12} \rho_{n}^{}\,A_{\text{w}}^{n}(\,t\,).
\end{equation}
Here, $A_{\text{w}}^{}$ denotes the effective wave amplitude. $\rho_{i}^{}\,(i=1,2,\cdots,12)$ denotes the coefficient of the i-th component of the polynomial that fits the relationship between $\Delta$GM and wave height at amidships in Fig.~\ref{fig:HBL_vs_GM}. The restoring arm for a ship heeling by two degrees in a regular wave is calculated based on the Froude-Krylov hypothesis~\cite{Hamamoto1991,Umeda1992} and a wavelength equal to the ship’s length is used. The GM for each wave amplitude is then calculated. The wave amplitude is positive when the wave trough is located at amidships, and the wave amplitude is negative when the wave crest is located at amidships, as shown in Fig.~\ref{fig:HBL_vs_GM}.

Furthermore, the methodology to derive the moment values of the roll angular acceleration is explained. The roll angular acceleration can be divided into two components as follows:
    \begin{equation}
        \begin{split}
        \label{accel_comp_all_re}
            \displaystyle
            &\ddot{\phi} = K_{1}^{}(\,\phi,\dot{\phi}\,) + K_{2}^{}(\,\phi,A_{\text{w}}^{}\,)
            \\
            \displaystyle
            \text{where} & \quad K_{1}^{}(\,\phi,\dot{\phi}\,) = - \beta_1^{}\,\dot{\phi} - \beta_3^{}\,\dot{\phi}_{}^3 - \sum_{n=1}^{5} \alpha_{2n-1}^{}\,\phi_{}^{2n-1}
            \\
            \displaystyle
            & \quad K_{2}^{}(\,\phi,A_{\text{w}}^{}\,) = - \frac{\omega_{0}^{2}}{GM} \sum_{n=1}^{12} \rho_{n}^{}\,A_{\text{w}}^{n}\, \phi
        \end{split} 
    \end{equation}
The roll angular acceleration is the sum of the stochastic variables $K_{1}^{}$ and $K_{2}^{}$. The first and second moments of the roll angular acceleration are represented by Eq.(\ref{accel_moment_first}) and (\ref{accel_moment_second}), respectively. The linearity of expectation can be used in this case.
    \begin{equation}
    \label{accel_moment_first}
        \mathbb{E}\left[\,\ddot{\phi}\,\right] = \mathbb{E}\left[\,K_1^{} + K_2^{}\,\right] = \mathbb{E}\left[\,K_1^{}\,\right] + \mathbb{E}\left[\,K_2^{}\,\right]
    \end{equation}
    \begin{equation}
    \label{accel_moment_second}
        \begin{split}
            \displaystyle
            \mathbb{E}\left[\,\ddot{\phi}_{}^{2}\,\right] & = \mathbb{E}\left[\,\left(\,K_1^{} + K_2^{}\,\right)_{}^{2}\,\right]
            \\
            \displaystyle
            & = \mathbb{E}\left[\,K_1^{2} + 2\,K_1^{}\,K_2^{} + K_2^{2} \,\right]
            \\
            \displaystyle
            & = \mathbb{E}\left[\,K_1^{2}\,\right] + 2\,\mathbb{E}\left[\,K_1^{}\,K_2^{}\,\right] + \mathbb{E}\left[\,K_2^{2}\,\right]
        \end{split} 
    \end{equation}
    where
    \begin{equation}
    \label{K1_moment_first}
    \displaystyle
        \mathbb{E}\left[\,K_1^{}\,\right] = - \beta_1^{}\,\mathbb{E}\left[\,\dot{\phi}\,\right] - \beta_3^{}\,\mathbb{E}\left[\, \dot{\phi}_{}^3\,\right] - \sum_{n=1}^{5} \alpha_{2n-1}^{}\,\mathbb{E}\left[\,\phi_{}^{2n-1}\,\right]
    \end{equation}
    \begin{equation}
    \label{K2_moment_first}
    \displaystyle
        \mathbb{E}\left[\,K_2^{}\,\right] = - \frac{\omega_{0}^{2}}{GM} \sum_{n=1}^{12} \rho_{n}^{}\,\mathbb{E}\left[\,\phi \,A_{\text{w}}^{n}\,\right]
    \end{equation}
    \begin{equation}
    \label{K1_moment_second}
        \begin{split}
            \displaystyle
            \mathbb{E}\left[\,K_1^{2}\,\right] = & \beta_1^{2}\,\mathbb{E}\left[\,\dot{\phi}_{}^{2}\,\right] + 2\,\beta_1^{}\, \beta_3^{}\,\mathbb{E}\left[\,\dot{\phi}_{}^{4}\,\right] + \beta_3^{2}\,\mathbb{E}\left[\,\dot{\phi}_{}^{6}\,\right] 
            \\
            \displaystyle
            & + \sum_{i=1}^{5}\sum_{j=1}^{5} \alpha_{2i-1}^{}\,\alpha_{2j-1}^{}\,\mathbb{E}\left[\,\phi_{}^{2i-1}\,\phi_{}^{2j-1}\, \right]
        \end{split}
    \end{equation}
    \begin{equation}
    \label{K2_moment_second}
        \displaystyle
        \mathbb{E}\left[\,K_2^{2}\,\right] = \frac{\omega_{0}^{4}}{GM_{}^{2}} \sum_{i=1}^{12}\sum_{j=1}^{12} \rho_{i}^{}\,\rho_{j}^{}\,\mathbb{E}\left[\,\phi_{}^{2}\,A_{\text{w}}^{i}\,A_{\text{w}}^{j}\,\right]
    \end{equation}
    \begin{equation}
    \label{K1K2_moment}
        \begin{split}
            \displaystyle
            \mathbb{E}\left[\,K_1^{}\,K_2^{}\,\right] = & \frac{\omega_{0}^{2}}{GM} \sum_{i=1}^{5}\sum_{j=1}^{12} \alpha_{2i-1}^{}\rho_{j}^{}\,\mathbb{E}\left[\,\phi_{}^{2i}\,A_{\text{w}}^{j}\,\right]
            \\
            \displaystyle
            & + \beta_1^{}\,\frac{\omega_{0}^{2}}{GM} \sum_{n=1}^{12} \rho_{n}^{}\,\mathbb{E}\left[\,\phi\,\dot{\phi}\, A_{\text{w}}^{n}\,\right]
            \\
            \displaystyle
            & + \beta_3^{}\,\frac{\omega_{0}^{2}}{GM} \sum_{n=1}^{12} \rho_{n}^{}\,\mathbb{E}\left[\,\phi\,\dot{\phi}_{}^3 \,A_{\text{w}}^{n}\,\right]
        \end{split}
    \end{equation}
Here, $\mathbb{E}[\,\cdot\,]$ denotes the moment of random values. The moment of random variables $x_{1}^{}x_{2}^{}\cdots$ is defined as follows:
    \begin{equation}
        \begin{split}
        \displaystyle
            & \mathbb{E}\left[\,x_{1}^{k_1}\,x_{2}^{k_2}\,\cdots\,\right] 
            \\
            \displaystyle
            &
            = \int_{-\,\infty}^{+\infty}\cdots\int_{-\,\infty}^{+\infty} x_{1}^{k_1}\,x_{2}^{k_2}\,\cdots\, p(\,x_{1}^{},x_{2}^{},\cdots\,)\,dx_{1}^{}\,dx_{2}^{}\,\cdots  .
        \end{split}
    \end{equation}
In this study, higher-order moments must be truncated. Therefore, the third and higher-order moment values included in Eq.~(\ref{accel_moment_first}) - (\ref{K1K2_moment}) can be approximated using the 2nd-order cumulant neglect closure method (Gaussian closure method)\cite{Sun1987}\cite{Sun1989}\cite{Wojtkiewicz1996}. Here, the moment values of the roll angle, roll velocity and effective wave amplitude are required. These moments can be obtained using our previously proposed methodology~\cite{Maruyama2022moment}.

\begin{figure*}[h]
    \centering
    \includegraphics[scale=0.55]{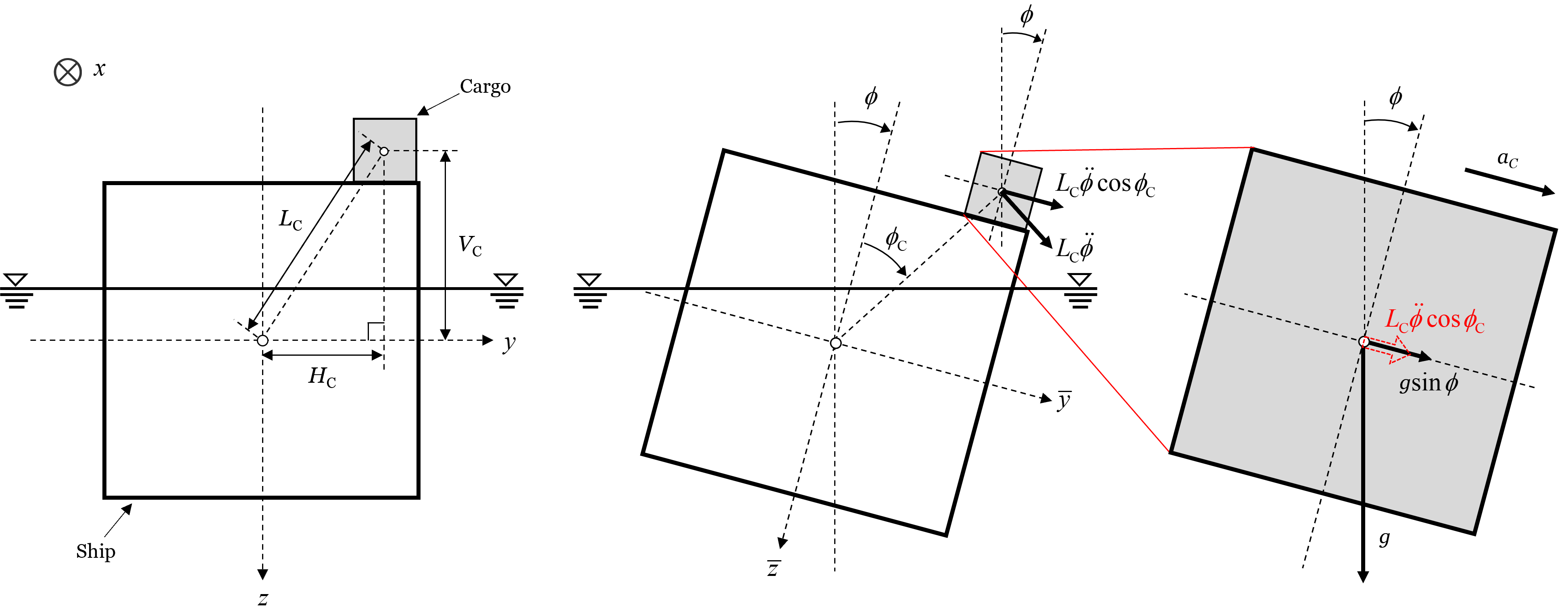}
    \caption{Coordinate system and schematic; this figure was adapted from Maruyama et al.~\cite{Maruyama2022accel}}
    \label{fig:cargo_accel}
\end{figure*}

\renewcommand{\arraystretch}{1.5}
\begin{table}[h]
    \begin{center}
    \caption{Positions of an arbitrary body on the subject ship}
        \begin{tabular}{llll}
        \hline
         {Position} & {$V_{\text{C}}^{}$ [\,m\,]} & {$H_{\text{C}}^{}$ [\,m\,]} & {$L_{\text{C}}^{}$ [\,m\,]}
         \\
         \hline
          $\text{C}_{1}^{}$ & $10.0$ & $0.00$ & $10.0$
          \\
          $\text{C}_{2}^{}$ & $30.0$ & $20.0$ & $36.1$
          \\
        \hline
        \end{tabular}
        \label{tab:Positions of an arbitrary body}
      \end{center}
\end{table}
\renewcommand{\arraystretch}{1.0}

Furthermore, the methodology for calculating the moment values of the cargo lateral acceleration is explained. First, the lateral acceleration $a_{\text{C}}^{}$ is represented in Eq.~(\ref{accel_cargo_all}). These relationship equations are derived by Maruyama et al.~\cite{Maruyama2022accel}. A cargo’s position is set in Table\ref{tab:Positions of an arbitrary body} in this study. The variables $\textstyle L_{\text{C}}^{}$, $\textstyle V_{\text{C}}^{}$, and $\textstyle H_{\text{C}}^{}$ are the direct, vertical, and horizontal distances, respectively.
\begin{equation}
    \label{accel_cargo_all}
        \displaystyle
        a_{\text{C}}^{} = K_{\text{C}1}(\,\phi,\dot{\phi}\,) + K_{\text{C}2}(\,\phi,A_{\text{w}}^{}\,)
\end{equation}
where
\begin{equation}
    \begin{split}
    \label{accel_cargo_all where}
        \displaystyle
        & K_{\text{C}1}(\,\phi,\dot{\phi}\,) = g\,\sin \phi + L'\,K_{1}(\,\phi,\dot{\phi}\,)
        \\
        \displaystyle
        & K_{\text{C}2}(\,\phi,A_{\text{w}}^{}\,) = L'\,K_{2}(\,\phi,A_{\text{w}}^{}\,)
        \\
        \displaystyle
        & L'= L_{\text{C}}^{}\,\cos{\phi_{\text{C}}^{}}
        \\
        \displaystyle
        & L_{\text{C}}^{}= \sqrt{V_{\text{C}}^{2} + H_{\text{C}}^{2}}
    \end{split}
\end{equation}
Here, $g$ denotes the gravitational acceleration. A cargo’s lateral acceleration is the sum of the stochastic variables $K_{\text{C}1}^{}$ and $K_{\text{C}2}^{}$. Therefore, the linearity of expectation can be used. The first and second moments of the cargo lateral acceleration is represented in Eqs.~(\ref{accel_arbitrary_moment_first}) and (\ref{accel_arbitrary_moment_second}).
    \begin{equation}
    \label{accel_arbitrary_moment_first}
        \begin{split}
            \displaystyle
            \mathbb{E}\left[\,a_{\text{C}}^{}\,\right] & = \mathbb{E}\left[\,g\,\sin \phi + L'\,K_{1}^{} + L'\,K_{2}^{}\,\right]
            \\
            \displaystyle
            & = g\,\mathbb{E}\left[\,\sin \phi\,\right] + L'\,\mathbb{E}\left[\,K_{1}^{}\,\right] + L'\,\mathbb{E}\left[\,K_{2}^{} \,\right]
        \end{split} 
    \end{equation}
    \begin{equation}
        \begin{split}
        \label{accel_arbitrary_moment_second}
            \displaystyle
            \mathbb{E}\left[\,a_{\text{C}}^{2}\,\right] & = \mathbb{E}\left[\,\left(\,g\,\sin \phi + L'\,K_{1}^{} + L'\,K_{2}^{} \,\right)_{}^{2}\,\right]
            \\
            \displaystyle
            &
            = L^{'2}\,\mathbb{E}\left[\,K_{1}^{2}\,\right] + 2\,L^{'2}\,\mathbb{E}\left[\,K_{1}^{}\,K_{2}^{}\,\right] + L^{'2}\, \mathbb{E}\left[\,K_{2}^{2}\,\right] 
            \\
            \displaystyle
            &
            + 2\,g\,L'\,\mathbb{E}\left[\,K_{1}^{}\,\sin{\phi}\,\right] + 2\,g\,L'\,\mathbb{E}\left[\,K_{2}^{}\,\sin{\phi}\,\right]
            \\
            \displaystyle
            &
            + g^{2}\,\mathbb{E}\left[\,\sin^{2} {\phi}\,\right]
        \end{split} 
    \end{equation}
In Eqs.~(\ref{accel_arbitrary_moment_first}) and (\ref{accel_arbitrary_moment_second}), the approach to the moment that includes the sinusoidal function must be considered. In this study, the series expansion of the sinusoidal function is used.
    \begin{equation}
    \label{series expansion of sinusoidal function}
        \sin{\phi} = \phi - \frac{\phi_{}^{3}}{6} + \frac{\phi_{}^{5}}{120} - \frac{\phi_{}^{7}}{5040} + \frac{\phi_{}^{9}}{362880} + O [\,\phi\,]_{}^{11} \,\, .
    \end{equation}
The extent of the remaining terms of the series expansion should be determined. Thus, the appropriate number of the term is examined by comparing the moment values obtained by solving Eqs.~(\ref{accel_arbitrary_moment_first}) and (\ref{accel_arbitrary_moment_second}). Here, the MCS result is used to calculate the moment values of the roll angle, roll velocity and effective wave amplitude. Therefore, Table\ref{tab:Moment values at a cargo of C1 checkorder} shows that the moment values converge sufficiently when the third or higher term in Eq.~(\ref{series expansion of sinusoidal function}) is used.
    \renewcommand{\arraystretch}{1.5}
    \begin{table}[h]
        \begin{center}
            \caption{Moment values at a cargo of Position:$\text{C}_{1}^{}$}
            \begin{tabular}{llll}
            \hline
             {order} & {$\mathbb{E}\left[K_{\text{C1}}^{}\right]$} & {$\mathbb{E}\left[K_{\text{C2}}^{}\right]$} & {$\mathbb{E}\left[a_{\text{C}}^{}\right]$} \\
            \hline
              $O[\,\phi_{}^{3}\,]$ & $-3.653\times10_{}^{-4}$ & $-1.736\times10_{}^{-6}$ & $-3.670\times10_{}^{-4}$ \\
              $O[\,\phi_{}^{5}\,]$ & $-3.567\times10_{}^{-4}$ & $-1.736\times10_{}^{-6}$ & $-3.584\times10_{}^{-4}$ \\
              $O[\,\phi_{}^{7}\,]$ & $-3.568\times10_{}^{-4}$ & $-1.736\times10_{}^{-6}$ & $-3.585\times10_{}^{-4}$ \\
              $O[\,\phi_{}^{9}\,]$ & $-3.568\times10_{}^{-4}$ & $-1.736\times10_{}^{-6}$ & $-3.585\times10_{}^{-4}$ \\
              $O[\,\phi_{}^{11}\,]$ & $-3.568\times10_{}^{-4}$ & $-1.736\times10_{}^{-6}$ & $-3.585\times10_{}^{-4}$ \\
            \hline
            \\
            \end{tabular}
            \begin{tabular}{llll}
            \hline
             {order} & {$\mathbb{E}\left[\,K_{\text{C1}}^{2}\,\right]$} & {$\mathbb{E}\left[\,K_{\text{C2}}^{2}\,\right]$} & {$\mathbb{E}\left[\,a_{\text{C}}^{2}\,\right]$} \\
            \hline
              $O[\,\phi_{}^{3}\,]$ & $3.683$ & $1.684\times10_{}^{-3}$ & $3.675$ \\
              $O[\,\phi_{}^{5}\,]$ & $3.514$ & $1.684\times10_{}^{-3}$ & $3.507$ \\
              $O[\,\phi_{}^{7}\,]$ & $3.516$ & $1.684\times10_{}^{-3}$ & $3.508$ \\
              $O[\,\phi_{}^{9}\,]$ & $3.516$ & $1.684\times10_{}^{-3}$ & $3.508$ \\
              $O[\,\phi_{}^{11}\,]$ & $3.516$ & $1.684\times10_{}^{-3}$ & $3.508$ \\
            \hline
            \end{tabular}
            \label{tab:Moment values at a cargo of C1 checkorder}
        \end{center}
    \end{table}
    \renewcommand{\arraystretch}{1.0}

%
\section{Calculation Results}
\subsection{Moment values of acceleration}\label{sec:sec4sub1}

When calculating the moment values of the roll angular and lateral accelerations, the moment values of the roll angle, roll velocity, and effective wave are required. The moment values in this study are of two types. One set of moment values is obtained from the MCS result, whereas the other is obtained from solving the moment equations. The latter has been conducted by Maruyama et al.~\cite{Maruyama2022moment}. These values are presented in Table\ref{tab:calculation_results_superposition_and_ME}. 
In Table\ref{tab:accel_superposition_and_moment_value_C11}, the result obtained by substituting the values of Table\ref{tab:calculation_results_superposition_and_ME} into Eqs.~(\ref{accel_moment_first}) and (\ref{accel_moment_second}) are presented. In this calculation, the second order cumulant neglect closure method is used. First, ``MCS" indicates the moment values derived from the MCS result directly. Second, ``Based MCS" indicates the moment values obtained by solving Eqs.~(\ref{accel_moment_first}) and (\ref{accel_moment_second}), based on the moment values of ``MCS" in Table\ref{tab:calculation_results_superposition_and_ME}. Finally, ``Moment eq.(2nd)" indicates the moment values obtained by solving Eqs.(\ref{accel_moment_first}) and (\ref{accel_moment_second}), based on the moment values of ``Moment eq.(2nd)" in Table\ref{tab:calculation_results_superposition_and_ME}.
The values of ``MCS" and ``Based MCS" are not much different, as shown by the second order moment presented in Table\ref{tab:accel_superposition_and_moment_value_C11}. However, the values of ``MCS" and ``Moment eq.(2nd)" have a discrepancy because the second order moment value of the roll angle and roll velocity differs between ``MCS" and ``Moment eq.(2nd)" in Table\ref{tab:calculation_results_superposition_and_ME}. Furthermore, a discrepancy exists between ``MCS", ``Based MCS", and ``Moment eq.(2nd)" for the fourth-order moment. The non-Gaussian does not affect these moment values when using the second order cumulant closure method. Moreover, the moment values of the roll angle and roll velocity derive from the moment equations are not shown in Table\ref{tab:calculation_results_superposition_and_ME}, which discrepancy in the MCS result. The first- and third-order moments are almost zero, whereas the second- and fourth-order moments are important while considering the PDF shape even if their magnitude is small. Therefore, in the future, we should conduct a study to reduce the discrepancy.

The moment of the cargo lateral acceleration is presented in Table\ref{tab:accel_arbitrary_superposition_and_moment_value_C11}. Compared among the three types of data, the same moment of the previously mentioned roll angular acceleration is obtained. ``Moment eq.(2nd)" in Table\ref{tab:accel_arbitrary_superposition_and_moment_value_C11} shows that the moment values of $\mathbb{E}\left[K_{\text{C1}}^{2}\right]$ and $\mathbb{E}\left[a_{\text{C}}^{2}\right]$ differ slightly, and this relationship is consistent with ``MCS". In other words, it can be observed that this relationship can be derived using moment equations.

    \renewcommand{\arraystretch}{1.5}
     \begin{table}[h]
      \begin{center}
        \caption{Moment values obtained using MCS and solving the moment equations. }
        \begin{tabular}{lll}
        \hline
         {} &  {MCS} &  {Moment eq.(2nd)}
         \\
         \hline
          $\mathbb{E}\left[\,\phi\,\right]$ & $-4.01\times10_{}^{-5}$ & $-6.22\times10_{}^{-4}$
          \\
          $\mathbb{E}\left[\,\phi_{}^{2}\,\right]$ & $4.38\times10_{}^{-2}$ & $3.62\times10_{}^{-2}$
          \\
          $\mathbb{E}\left[\,\dot{\phi}_{}^{}\,\right]$ & $-6.00\times10_{}^{-6}$ & $-2.25\times10_{}^{-4}$
          \\
          $\mathbb{E}\left[\,\dot{\phi}_{}^{2}\,\right]$ & $2.81\times10_{}^{-3}$ & $2.32\times10_{}^{-3}$
          \\
          $\mathbb{E}\left[\,A_{\text{w}}^{}\,\right]$ & $-1.66\times10_{}^{-5}$ & $0.00$
          \\
          $\mathbb{E}\left[\,A_{\text{w}}^{2}\,\right]$ & $0.786$ & $0.794$
          \\
          $\mathbb{E}\left[\,\phi\,\dot{\phi}_{}^{}\,\right]$ & $-7.32\times10_{}^{-7}$ & $1.56\times10_{}^{-5}$
          \\
          $\mathbb{E}\left[\,\phi\,A_{\text{w}}^{}\,\right]$ & $1.01\times10_{}^{-5}$ & $2.54\times10_{}^{-4}$
          \\
          $\mathbb{E}\left[\,\dot{\phi}_{}^{}\,A_{\text{w}}^{}\,\right]$ & $-8.33\times10_{}^{-7}$ & $3.04\times10_{}^{-5}$
          \\
        \hline
        \end{tabular}
        \label{tab:calculation_results_superposition_and_ME}
      \end{center}
    \end{table}
    \renewcommand{\arraystretch}{1.0}
    
    \renewcommand{\arraystretch}{1.5}
    \begin{table}[h]
        \begin{center}
            \caption{MCS result obtained using the superposition principle and calculation results from solving Eqs.~(\ref{accel_moment_first}) and (\ref{accel_moment_second}) }
            \begin{tabular}{llll}
            \hline
            {} & { MCS } & { Based MCS } & { Moment eq.(2nd) } \\
            \hline
              $\mathbb{E}\left[\,\ddot{\phi}_{}^{}\,\right]$ & $1.62\times10_{}^{-6}$ & $2.63\times10_{}^{-6}$ & $3.76\times10_{}^{-6}$ \\
              $\mathbb{E}\left[\,\ddot{\phi}_{}^{2}\,\right]$ & $2.07\times10_{}^{-4}$ & $2.09\times10_{}^{-4}$ & $1.71\times10_{}^{-4}$ \\
              $\mathbb{E}\left[\,\ddot{\phi}_{}^{3}\,\right]$ & $-2.58\times10_{}^{-9}$ & $1.01\times10_{}^{-9}$ & $2.17\times10_{}^{-8}$ \\
              $\mathbb{E}\left[\,\ddot{\phi}_{}^{4}\,\right]$ & $1.69\times10_{}^{-7}$ & $1.29\times10_{}^{-7}$ & $8.73\times10_{}^{-8}$ \\
            \hline
            \end{tabular}
            \label{tab:accel_superposition_and_moment_value_C11}
        \end{center}
    \end{table}
    \renewcommand{\arraystretch}{1.0}
    
    \renewcommand{\arraystretch}{1.5}
    \begin{table}[h]
        \begin{center}
            \caption{MCS result obtained using the superposition principle and calculation results from solving Eqs.~(\ref{accel_arbitrary_moment_first}) and (\ref{accel_arbitrary_moment_second})}
            \begin{tabular}{llll}
            {} & {} & {} & {Position : C1} \\
            \hline
            {} & { MCS } & { Based MCS } & { Moment eq.(2nd) } \\
            \hline
              $\mathbb{E}\left[\,K_{\text{C1}}^{}\,\right]$ & $-1.79\times10_{}^{-4}$ & $-3.57\times10_{}^{-4}$ & $-5.52\times10_{}^{-3}$ \\
              $\mathbb{E}\left[\,K_{\text{C2}}^{}\,\right]$ & $-2.85\times10_{}^{-5}$ & $-1.74\times10_{}^{-6}$ & $-4.83\times10_{}^{-5}$ \\
              $\mathbb{E}\left[\,a_{\text{C}}^{}\,\right]$ & $-1.79\times10_{}^{-4}$ & $-3.58\times10_{}^{-4}$ & $-5.57\times10_{}^{-3}$ \\
              $\mathbb{E}\left[\,K_{\text{C1}}^{2}\,\right]$ & $3.69$ & $3.52$ & $2.93$ \\
              $\mathbb{E}\left[\,K_{\text{C2}}^{2}\,\right]$ & $3.74\times10_{}^{-3}$ & $1.68\times10_{}^{-3}$ & $1.40\times10_{}^{-3}$ \\
              $\mathbb{E}\left[\,a_{\text{C}}^{2}\,\right]$ & $3.69$ & $3.51$ & $2.92$ \\
            \hline
            \\
            {} & {} & {} & {Position : C2} \\
            \hline
            {} & { MCS } & { Based MCS } & { Moment eq.(2nd) } \\
            \hline
              $\mathbb{E}\left[\,K_{\text{C1}}^{}\,\right]$ & $-1.36\times10_{}^{-4}$ & $-3.09\times10_{}^{-4}$ & $-4.58\times10_{}^{-3}$ \\
              $\mathbb{E}\left[\,K_{\text{C2}}^{}\,\right]$ & $-2.01\times10_{}^{-5}$ & $-5.51\times10_{}^{-6}$ & $-1.45\times10_{}^{-4}$ \\
              $\mathbb{E}\left[\,a_{\text{C}}^{}\,\right]$ & $-1.42\times10_{}^{-4}$ & $-3.14\times10_{}^{-4}$ & $-4.72\times10_{}^{-3}$ \\
              $\mathbb{E}\left[\,K_{\text{C1}}^{2}\,\right]$ & $2.73$ & $2.73$ & $2.16$ \\
              $\mathbb{E}\left[\,K_{\text{C2}}^{2}\,\right]$ & $2.00\times10_{}^{-2}$ & $1.52\times10_{}^{-2}$ & $1.26\times10_{}^{-2}$ \\
              $\mathbb{E}\left[\,a_{\text{C}}^{2}\,\right]$ & $2.73$ & $2.73$ & $2.15$ \\
            \hline
            \end{tabular}
            \label{tab:accel_arbitrary_superposition_and_moment_value_C11}
        \end{center}
    \end{table}
    \renewcommand{\arraystretch}{1.0}

\subsection{Procedure for determining the PDF, and results}\label{sec:sec4sub2}
The PDF of roll angular acceleration or cargo lateral acceleration is calculated based on the moment values presented in Tables\ref{tab:accel_superposition_and_moment_value_C11} and \ref{tab:accel_arbitrary_superposition_and_moment_value_C11}. The following non-Gaussian PDF shape types are set in this study. The shape of type1 is inspired from the Laplace distribution, and type2 is inspired from a logistic distribution.
    \begin{equation}
        \label{PDF_opt_type1}
        \begin{split}
            \displaystyle
            &\text{type1 : }
            \\
            \displaystyle
            & \mathcal{P}_{}^{}(\,X\,) = C \exp{ \left\{ - \left( d_1^{}\left| \,X\, \right| + d_2^{} \left| \,X\, \right|_{}^{2} + d_3^{} \left| \,X\, \right|_{}^{3} + d_4^{} \left| \,X\, \right|_{}^{4} \right) \right\} }
        \end{split}
    \end{equation}
    \begin{equation}
        \label{PDF_opt_type2}
        \begin{split}
            \displaystyle
            &\text{type2 : }
            \\
            \displaystyle
            & \mathcal{P}_{}^{}(\,X\,) = C \dfrac{ \exp{ \left( - \dfrac{ \,X\, }{ d_1^{} } \right) } }{ \left( 1 + \exp{ \left( - \dfrac{ \,X\, }{ d_1^{} } \right) } \right)_{}^{2} }
        \end{split}
    \end{equation}
    Here, $C$ denotes a normalization constant. In this study, we assumed that the first and third-order moments are almost zero, whereas the second and fourth-order moments are important for deriving variance and kurtosis. Therefore, the following expression is proposed to determine the coefficients of Eqs.~(\ref{PDF_opt_type1}) and (\ref{PDF_opt_type2}).
    \begin{equation}
        \label{Estimation_PDF_arbitrary1}
        \displaystyle
        J_{n}^{} = \left\{ 
        \begin{split}
            \displaystyle
            & \int_{-\,\infty}^{+\infty} \ddot{X}_{}^{n} \mathcal{P}(\,X\,)\,\mathrm{d}X - \mathbb{E}\left[ \,\ddot{X}_{}^{n}\, \right] \quad \left( n \text{ : odd }  \right)
            \\\\
            \displaystyle
            & \dfrac{ \displaystyle \int_{-\,\infty}^{+\infty} \ddot{X}_{}^{n} \mathcal{P}(\,X\,)\,\mathrm{d}X - \mathbb{E}\left[ \,\ddot{X}_{}^{n}\, \right]}{\mathbb{E}\left[\, \ddot{X}_{}^{n}\, \right]} \quad \left( n \text{ : even }  \right)
        \end{split}
        \right.
    \end{equation}
    The moment values obtained from Eqs.~(\ref{accel_moment_first}) and (\ref{accel_moment_second}) or Eqs.~(\ref{accel_arbitrary_moment_first}) and (\ref{accel_arbitrary_moment_second}) are only the first and second order moments. Thereby, the cumulant neglect closure method is used to obtain the higher-order moment values. Furthermore, the following objective function $\textstyle J\left( d_{1}^{},d_{2}^{},d_{3}^{},d_{4}^{}\right)$ is set.
    \begin{equation}
        \label{Estimation_PDF_arbitrary2}
        J\left(\,d_{1}^{},d_{2}^{},d_{3}^{},d_{4}^{}\,\right) = \sum_{i=1}^{4} l_{i}^{}\,|\,J_{i}^{}\,|
    \end{equation}
    Here, $\textstyle l_{i}^{}$ are weights. The values of type1 and type2 are $\textstyle l_{i}^{}=1\,(i = 1, 2, 3, 4)$ and $\textstyle l_{1}^{}=1,\, l_{i}^{}=0\,(i = 2, 3, 4)$, respectively.

    \begin{figure*}[tb]
        \centering
        \includegraphics[scale=0.8]{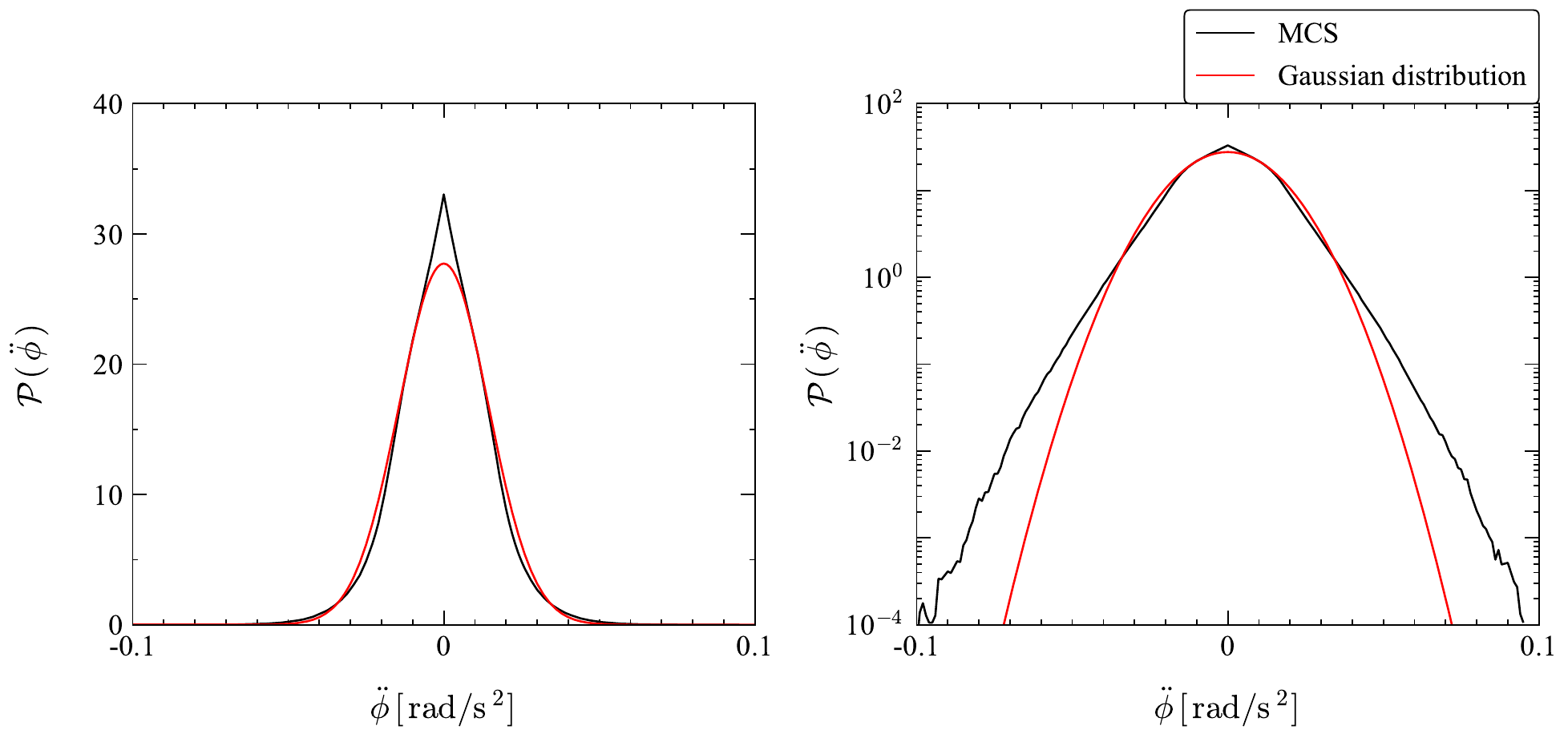}
        \caption{Comparing the PDF of roll angular acceleration between the MCS result and Gaussian distribution.}
        \label{fig:PDF_of_opt4}
    \end{figure*}
    \begin{figure*}[tb]
        \centering
        \includegraphics[scale=0.8]{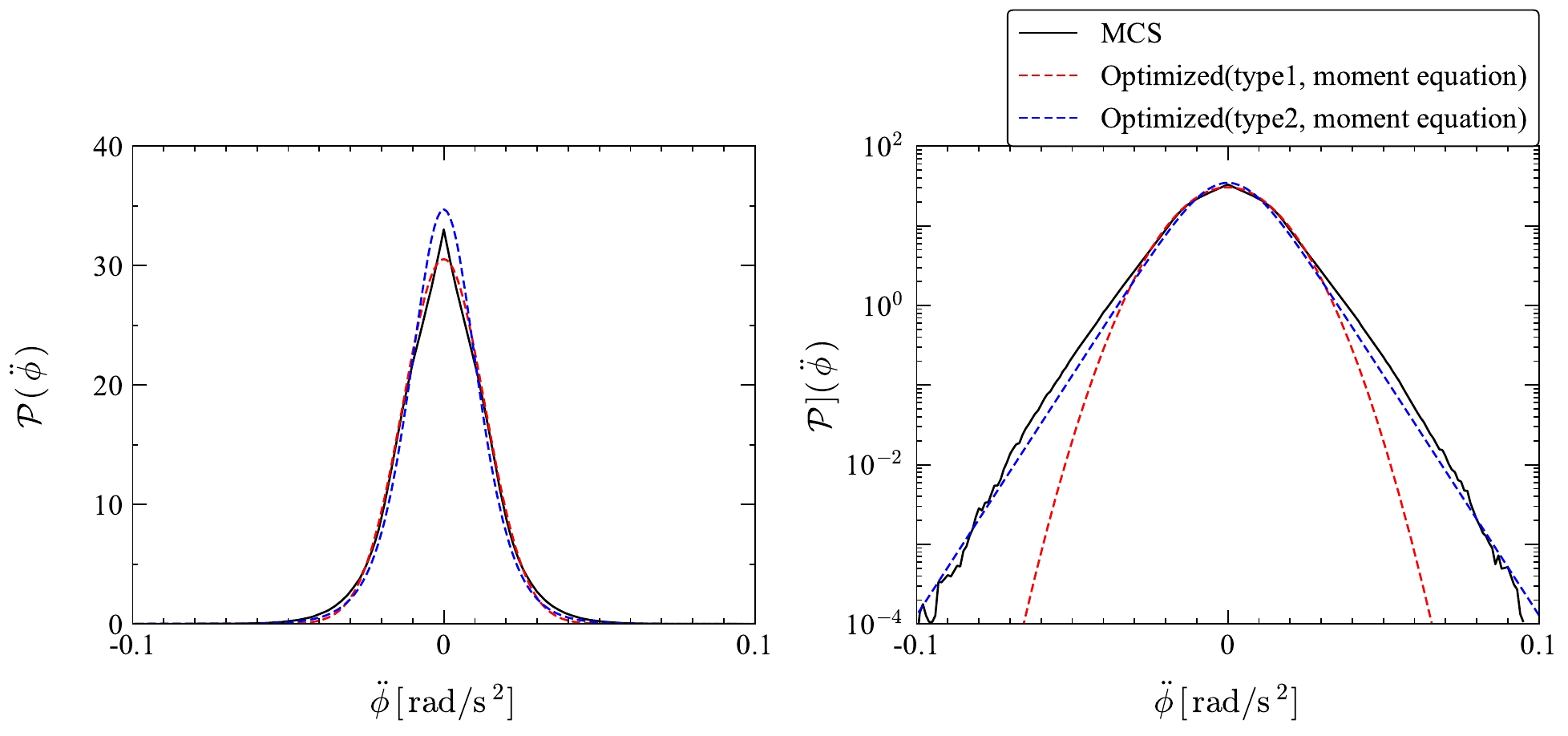}
        \caption{Comparing the PDF of roll angular acceleration between the MCS and optimized results using Eqs.~(\ref{PDF_opt_type1}) and (\ref{PDF_opt_type2}). Here, the moment values (Moment eq.(2nd)) presented in Table \ref{tab:accel_superposition_and_moment_value_C11} are used.}
        \label{fig:PDF_of_opt1}
    \end{figure*}
    \begin{figure*}[tb]
        \centering
        \includegraphics[scale=0.8]{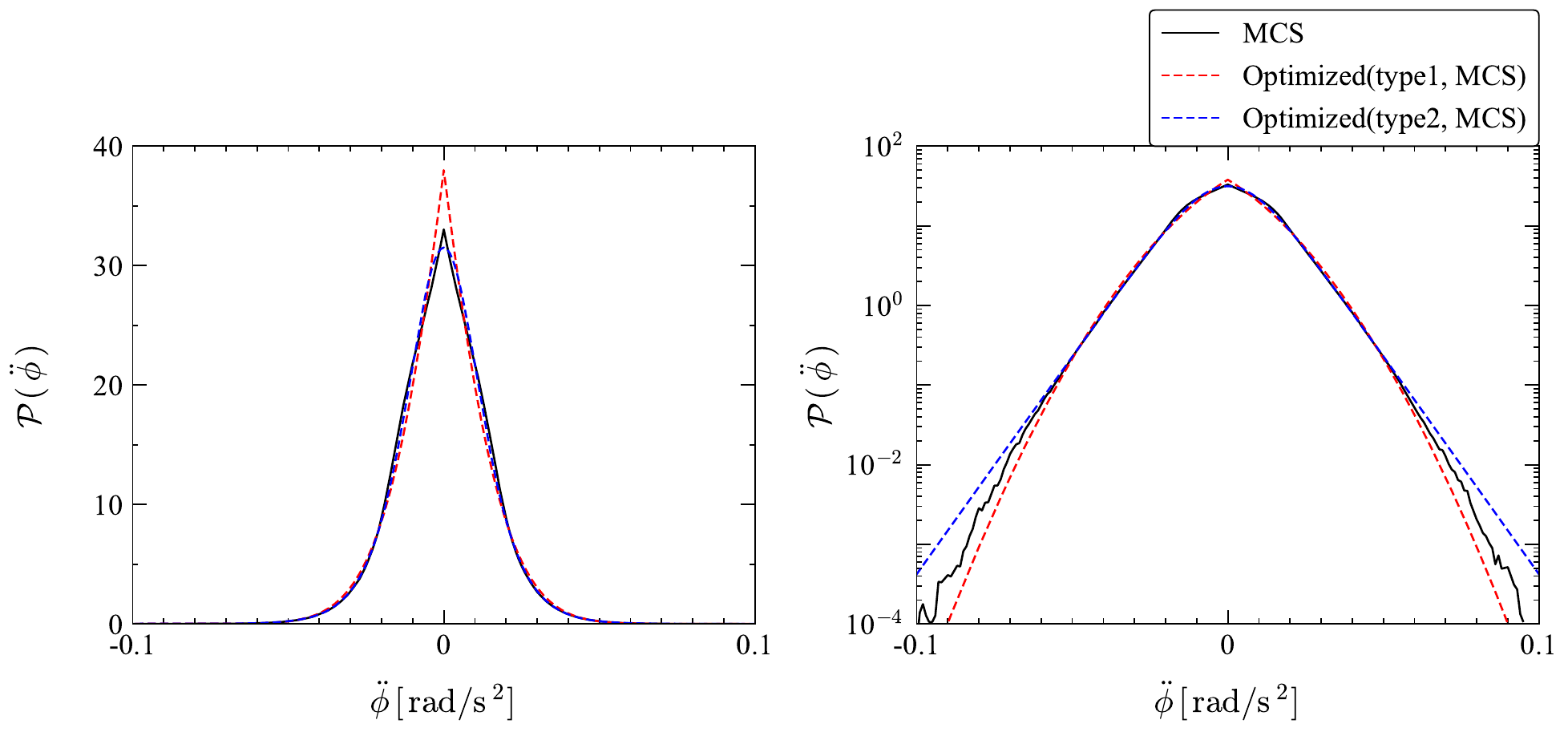}
        \caption{Comparing the PDF of roll angular acceleration between the MCS and optimized results using Eqs.~(\ref{PDF_opt_type1}) and (\ref{PDF_opt_type2}). Here, the moment values (MCS) presented in Table \ref{tab:accel_superposition_and_moment_value_C11} are used.}
        \label{fig:PDF_of_opt2}
    \end{figure*}
    \begin{figure*}[tb]
        \centering
        \includegraphics[scale=0.8]{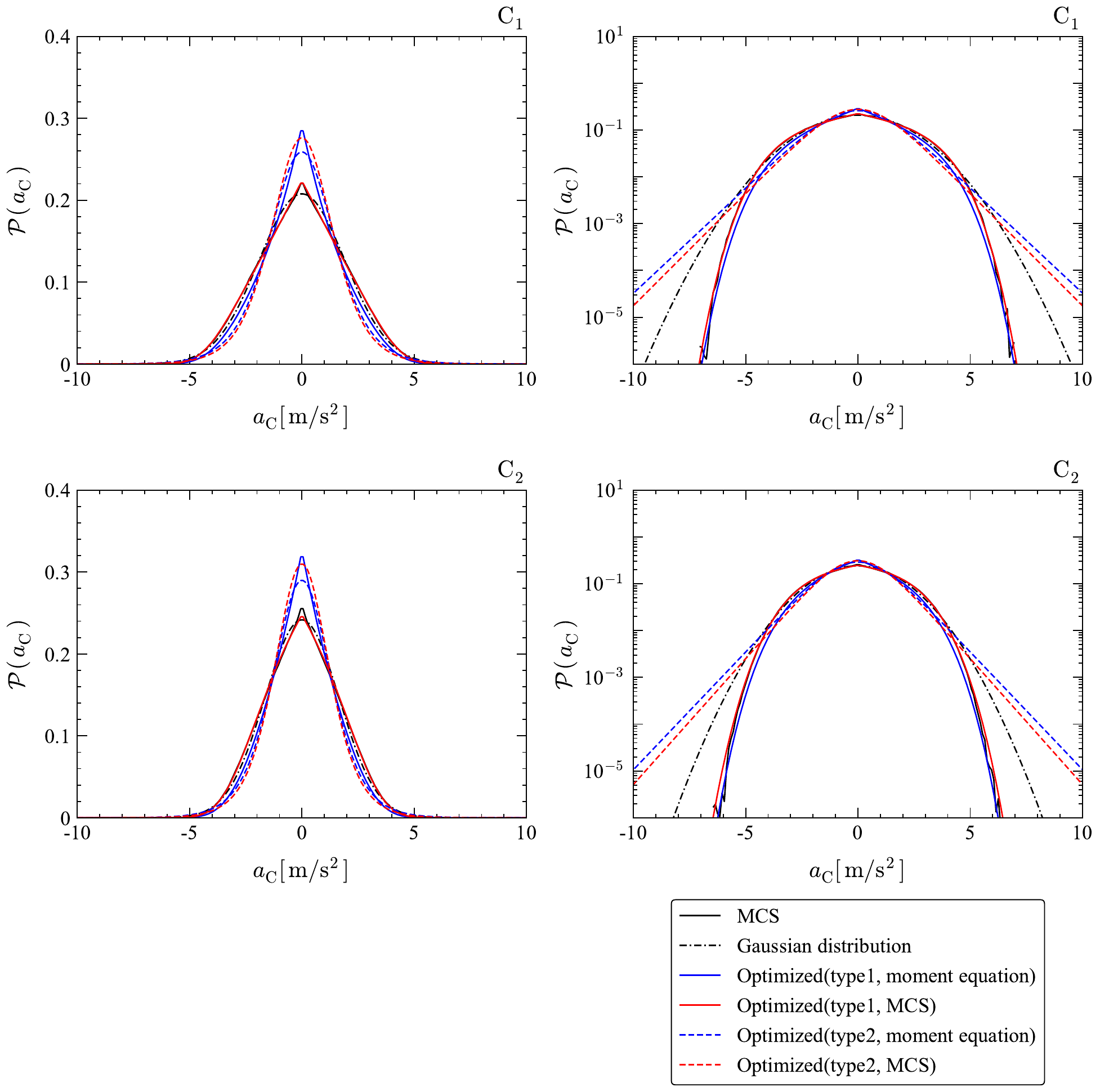}
        \caption{Comparing the PDF of a cargo’s lateral acceleration between the MCS result, Gaussian distribution, and optimized results using Eq.~(\ref{PDF_opt_type1}). Here, the moment values ( Moment eq.(2nd) or MCS ) presented in Table\ref{tab:accel_superposition_and_moment_value_C11} are used.}
        \label{fig:PDF_of_arbitrary_opt}
    \end{figure*}



Fig.\ref{fig:PDF_of_opt4} shows that the Gaussian distribution should not be used as the PDF for roll angular acceleration. Here, the mean and variance use the values obtained from the MCS result. The PDF of roll angular acceleration is determined using the moment values of ``Moment eq.(2nd)" presented in Table\ref{tab:accel_superposition_and_moment_value_C11}, as shown in Fig.\ref{fig:PDF_of_opt1}. In the range where the probability density low, type1 has a discrepancy with the MCS result. However, type2 has a close shape with the MCS result in all ranges.
To estimate the validity of the proposed PDF shape in this study, the PDF shape is determined based on the moment values of ``MCS" presented in Table\ref{tab:accel_superposition_and_moment_value_C11}. Both types have a close shape with the MCS result, as shown in Fig.\ref{fig:PDF_of_opt2}. When the PDF of roll angular acceleration is derived, type2 is a simpler and better PDF shape than type1. 

Furthermore, Fig.\ref{fig:PDF_of_arbitrary_opt} shows the PDF of roll angular acceleration determined using the moment values of ``Moment eq.(2nd)" and ``MCS" presented in Table\ref{tab:accel_arbitrary_superposition_and_moment_value_C11}. In the case of the Gaussian distribution, the values obtained from the MCS result are used for the mean and variance. When the PDF coefficient is determined based on the moment values of ``Moment eq.(2nd)," the PDF does not agree with the MCS result. Furthermore, the PDF generally agrees in the range where the probability density is low. Furthermore, the red solid line in Fig.~\ref{fig:PDF_of_arbitrary_opt} represents the optimized result when the coefficient of the PDF is determined using the moment values of ``MCS." Therefore, when the appropriate moment values and PDF shape of type1 are used, the PDF of a cargo lateral acceleration agreeing with the MCS result can be obtained. However, the PDF of type2 does not agree with the MCS result.
The above consideration is common in the two positions presented in Table\ref{tab:Positions of an arbitrary body}.  

The moment related to kurtosis is important as mentioned above. Thus, we should improve the calculation accuracy to obtain the appropriate moment values. Moreover, our proposed PDF shapes outperform those of the Gaussian PDF.

%
\section{Concluding Remarks}

Based on the linearity of expectation, the moment values of roll angular and cargo lateral accelerations can be derived using the moment values of the roll angle, roll velocity, and effective wave amplitude. If the appropriate moment values are used, the appropriate acceleration of the moment values can be obtained. Furthermore, we proposed a new PDF shape for the roll angular acceleration, which exhibits good agreement with the MCS result using the appropriate moment values. However, more research on the PDF shape and methodology moment values is required. We will derive higher-order moment values using the higher-order cumulant neglect closure method in the future.

%
\begin{acknowledgements}
This work was supported by a Grant-in-Aid for Scientific Research from the Japan Society for Promotion of Science (JSPS KAKENHI Grant Number 19H02360) and by JST SPRING, Grant Number JPMJSP2138, as well as the collaborative research program and financial support from the Japan Society of Naval Architects and Ocean Engineers.
This study was supported by the Fundamental Research Developing Association for Shipbuilding and Offshore (REDAS), managed by the Shipbuilders’ Association of Japan from April 2020 to March 2023. The authors would like to thank Enago(www.enago.jp) for English language review.
\end{acknowledgements}

\bibliographystyle{spphys}       
\bibliography{ref.bib}

\end{document}